\documentclass[a4wide,11pt]{scrartcl}

\usepackage[USenglish]{babel}
\usepackage[T1]{fontenc}
\usepackage[utf8]{inputenc}

\usepackage{graphicx,subfigure,tikz,pgfplots,epstopdf}
\graphicspath{{figs/}}

\usepackage{amsmath,amssymb,amsfonts,amsthm,array,geometry,stmaryrd,booktabs,paralist,todonotes,mathtools,pdflscape,siunitx,placeins}

\usepackage[pdftex,colorlinks=true,linkcolor=blue,citecolor=green,urlcolor=blue,bookmarks]{hyperref}

\numberwithin{equation}{section}

\newcommand{\R}{\mathbb{R}}
\newcommand{\dd}{\textup{d}}

\newtheorem{thm}{Theorem}[section]
\newtheorem{lem}[thm]{Lemma}
\newtheorem{defi}[thm]{Definition}
\newtheorem{rem}[thm]{Remark}

\newenvironment{prf}{\textbf{Proof:}} {\hspace*{\fill} $\square$ \newline}

\usepackage{url}
\usepackage{dsfont}
\usepackage{amsmath}
\usepackage{nicefrac}
\usepackage{xspace}
\usepackage{framed}

\usepackage{cite}

\newcommand{\p}[2]{\ensuremath{\frac{\partial #1}{\partial #2 }}}

\begin{document}
    \title{Boundary feedback control for hyperbolic systems}
    \author{Michael Herty\footnotemark[1] \and  Ferdinand Thein\footnotemark[1] \footnotemark[2]}
    \date{\today}
    \maketitle
    \begin{abstract}
        We are interested in the feedback stabilization of general linear multi-dimensional first order hyperbolic systems in $\R^d$.
        Using a Lyapunov function with a suited weight function depending on the system under consideration we show stabilization in $L^2$ for the studied system using a suitable feedback control.
        Therefore the controllability of the studied system is related to the feasibility of an associated linear matrix inequality.
        We show the applicability discussing the barotropic Euler equations.
    \end{abstract}
    \renewcommand{\thefootnote}{\fnsymbol{footnote}}
    \footnotetext[1]{IGPM RWTH Aachen, Templergraben 55, D-52056 Aachen, Germany.\\
    \href{mailto:herty@igpm.rwth-aachen.de}{\textit{herty@igpm.rwth-aachen.de}}}
    \footnotetext[2]{present address: Johannes Gutenberg-Universit\"at, Staudingerweg 9, D-55128 Mainz, Germany.\\
    \href{mailto:fthein@uni-mainz.de}{\textit{fthein@uni-mainz.de}}}

    \renewcommand{\thefootnote}{\arabic{footnote}}

    \newpage
    %
    \section{Introduction}
    The stabilization of spatially one--dimensional systems of hyperbolic balance laws is a vivid subject attracting research interest in the mathematical as well as in the engineering community
    and we refer to the monographs \cite{Bastin2016,MR2302744,MR2655971,MR2412038} for further references. The mentioned references also provide a comprehensive overview on related controllability problems.
    A particular focus has been put on problems modeled by the barotropic Euler equations and the shallow water equations which in one space dimension form a 2$\times$2 hyperbolic system to model the temporal
    and spatial evolution of fluid flows including flows on networks.
    Analytical results concerning the boundary control of such systems have been studied in several articles, cf.\ \cite{G1,G2,G3,G4} for gas flow and for water flow we refer to \cite{W1,W2,W3,W4,W5,W6}.
    One key aspect in the analysis is the Lyapunov function which is introduced as a weighted $L^2$ (or $H^s$) norm and which allows to estimate deviations from steady states, see e.g. \cite{Bastin2016}.
    Under rather general {dissipative} conditions the exponential decay of the Lyapunov function has been established for various problem formulations and we exemplary refer to \cite{L1,L2,L4,L5}.
    For a study on comparisons to other stability concepts we mention \cite{L7}.
    Stability with respect to a higher $H^s$-norm $(s\geq 2)$ gives stability of the nonlinear system \cite{L5,Bastin2016}.
    Without aiming at completeness, we mention that recently the results have been extended to also deal with e.g.\ input-to-state stability \cite{MR2899713},
    numerical methods have been discussed in \cite{MR3956429,MR3031137,MR3648349} and for results concerning nonlocal hyperbolic partial differential equations (PDEs) see \cite{MR4172728}.
    \par
    However, to the best of our knowledge the presented results are limited to the spatially one--dimensional case.
    Based on an example in metal forming processes, see \cite{Bambach2022,Herty2023} we extended results to multi--dimensional hyperbolic balance laws.
    A specific system in two dimensions is discussed in \cite{Dia2013} where a control problem for the shallow water equations is studied.
    There the authors take advantage of the structure of the system and show that the energy is non-increasing upon imposing suited boundary conditions.
    Just very recently a preprint was published where the boundary stabilization for two dimensional systems is studied using a different Lyapunov function, see \cite{Yang2023}
    and a comparative study is given in \cite{Herty2024}.
    A further recent result on multi--dimensional hyperbolic scalar conservation laws was presented in \cite{Serre2022}.
    However, the goal of the mentioned paper is different from the one presented here.
    In \cite{Serre2022} the author discusses scalar conservation laws in $\R^d$ whereas we are concerned with systems in a bounded domain.
    Moreover, we consider deviations from a reference state and thus study the linearized system, whereas in \cite{Serre2022} a situation is studied when two fixed states at infinity are prescribed.
    It will be interesting and part of future research to study the relation between the works \cite{Dia2013,Yang2023,Serre2022} and the results presented in this work.
    \par
    Throughout this work we will fundamentally rely on the assumption that the studied system is symmetric hyperbolic, i.e., the Jacobian matrices are symmetric.
    Symmetry of the system, in particular in multi-dimensions, gives well posed problems as it was shown in the pioneering work  \cite{Friedrichs1954}.
    This is of major importance since it is not possible to expect strict hyperbolicity, another property of hyperbolic systems leading to well posedness, for relevant systems in multi-dimensions.
    In \cite{Godunov1961} it is shown that the Euler equations can be rewritten as a symmetric hyperbolic system.
    Later it was proven that all hyperbolic systems which are equipped with an additional conservation law for a quantity which is convex in terms of the conserved quantities can be put into symmetric form, see
    \cite{Friedrichs1971}.
    Thus many, if not to say all, relevant systems arising in applications can be symmetrized since they enjoy this additional conservation law representing, e.g., the conservation of energy or entropy.
    For some examples on symmetrizable systems we refer to \cite{Boillat1996}.
    Also, well--posedness results related to symmetric hyperbolic systems are obtained in several publications and without assuming completeness we exemplary refer \cite{Brenner1966} for results in $L^p$
    and \cite{Kato1975} for the linear and quasi linear case. Initial boundary value problems have been studied e.g. in \cite{Majda1975} or \cite{Peyser1975} in the case of constant coefficients.
    For a comprehensive presentation of results on initial boundary value problems for Friedrichs symmetrizable multi-dimensional systems we refer to \cite{BenzoniGavage2007}.
    Further insight on the properties of these large class of systems are presented in \cite{Dafermos2016,Ruggeri2021} where a brief historical review of this topic is given in \cite{Ruggeri2021}.
    There are basically three major points where we make use of the symmetry property. First, it ensures well-posedness and the existence of sufficiently smooth solutions as already noted before.
    Second, the symmetric Jacobians allow to use the product rule for a quadratic form which is needed throughout the proof of the main result.
    The third point where we make use of this property is the discussion of a \emph{linear matrix inequality (LMI)} which is essential for a needed estimate 
    and we exemplary refer to \cite{Boyd1994} for a broad introduction and further references. 
    The symmetric hyperbolicity of the system induces a favorable algebraic structure which was already investigated by Lax in \cite{Lax1958}.
    This work stimulated further research and a related conjecture was answered in \cite{Helton2007}. We further recommend \cite{Vinnikov2012} for a comprehensive overview and additional references.
    \par
    In this work we will study the boundary stabilization of multi--dimensional linear hyperbolic balance laws with variable coefficient matrices using an extension of the Lyapunov function introduced in \cite{Bastin2016}.
    The particular structure of a hyperbolic system is exploited to derive a condition on the feedback law, such that the Lyapunov function decays exponentially fast, see Section \ref{sec:defthm}
    and Section \ref{sec:prf} for the proof.
    The studied results will be applied the barotropic Euler equations in Section \ref{sec:app_hyp_sys}.
    \section{Stabilization of Multi-Dimensional Linear Hyperbolic Balance Laws}\label{sec:defthm}
    We are interested in the following initial boundary value problem (IBVP) for the given system of hyperbolic PDEs
    \begin{align}
        \begin{dcases}
            \p{}{t}\mathbf{w}(t,\mathbf{x}) + \sum_{k=1}^d\mathbf{A}^{(k)}(\mathbf{x})\p{}{x_k}\mathbf{w}(t,\mathbf{x}) + \mathbf{B}(\mathbf{x})\mathbf{w}(t,\mathbf{x}) &= 0,\;(t,\mathbf{x})\in[0,T)\times\Omega\\
            \mathbf{w}(0,\mathbf{x}) &= \mathbf{w}_0(\mathbf{x}),\;\mathbf{x}\in\Omega,\\
            \mathbf{w}(t,\mathbf{x}) &= \mathbf{w}_{BC}(t,\mathbf{x}),\;(t,\mathbf{x})\in[0,T)\times\partial\Omega
        \end{dcases}\label{eq:hyp_cons_sys}
    \end{align}
    Here $\mathbf{w}(t,\mathbf{x}) \equiv (w_1(t,\mathbf{x}),\dots,w_n(t,\mathbf{x}))^T$ is the vector of unknowns
    and $\Omega \subset \R^{d}$ a bounded domain with piecewise $C^1$ smooth boundary $\partial\Omega$.
    Moreover, $\mathbf{A}^{(k)}(\mathbf{x})$ and $\mathbf{B}(\mathbf{x})$ are sufficiently smooth and bounded $n\times n$ real matrices. The $\mathbf{A}^{(k)}(\mathbf{x})$ are in particular symmetric.
    Upon applying a suited variable transformation we may transform any symmetric hyperbolic system into the given form with the identity in front of the time derivative, cf. \cite{BenzoniGavage2007}.
    The assumption of symmetry is no major restriction since it includes all systems equipped with an additional conservation law, cf. \cite{Friedrichs1971,Dafermos2016}.
    This includes most systems relevant for applications, see \cite{Boillat1996,Ruggeri2021}.
    We define the matrix
    \[
      \mathbf{A}^\ast(\mathbf{x},\mathbf{\nu}) := \sum_{k=1}^d \nu_k\mathbf{A}^{(k)}(\mathbf{x})
    \]
    with $\mathbf{\nu} = (\nu_1,\dots,\nu_d)^T \in \mathbb{S}^{d-1}$ being a unit vector in $\R^d$ and by construction $\mathbf{A}^\ast(\mathbf{x},\mathbf{\nu})$ is also symmetric.
    Adopting the notion given in \cite{Lukacova2006,Helton2007} we refer to the matrix $\mathbf{A}^\ast$ as the \emph{pencil matrix}.
    Concerning the definition of hyperbolicity of system \eqref{eq:hyp_cons_sys} we follow \cite{Dafermos2016}.
    System \eqref{eq:hyp_cons_sys} is said to be hyperbolic if the matrix $\mathbf{A}^\ast(\mathbf{x},\mathbf{\nu})$ has $n$ real eigenvalues $\lambda_i = \lambda_i(\mathbf{x},\mathbf{\nu}),\, i = 1,\dots, n$,
    and $n$ corresponding linearly independent right eigenvectors $\mathbf{r}_i = \mathbf{r}_i(\mathbf{\mathbf{x},\nu}), i = 1,\dots, n$ for all $\mathbf{\nu} \in \mathbb{S}^{d-1}$.
    Note that by choosing
    \[
      \mathbf{\nu} = \mathbf{e}_k = (0,\dots,0,\underbrace{1}_{k},0,\dots,0)^T
    \]
    we have $\mathbf{A}^\ast(\mathbf{x},\mathbf{e}_k) = \mathbf{A}^{(k)}(\mathbf{x})$ and thus the individual Jacobians $\mathbf{A}^{(k)}(\mathbf{x})$ are also diagonalizable with real eigenvalues.
    Further following \cite{BenzoniGavage2007} there exists a well-conditioned orthogonal matrix $\mathbf{T}(\mathbf{x},\mathbf{\nu})$ such that
    \begin{align}
        \mathbf{\Lambda}^\ast(\mathbf{x},\mathbf{\nu}) = \mathbf{T}^T(\mathbf{x},\mathbf{\nu})\mathbf{A}^\ast(\mathbf{x},\mathbf{\nu})\mathbf{T}(\mathbf{x},\mathbf{\nu})\label{eq:diag_transform}
    \end{align}
    is a diagonal matrix containing the eigenvalues of the system \eqref{eq:hyp_cons_sys}.
    The afore mentioned properties are distinctive features of multi--dimensional hyperbolic PDEs and we highly recommend the references \cite{BenzoniGavage2007,Dafermos2016}, cited above, for further reading.
    A further delicate issue is the dependence of the eigenspaces and eigenvalues on $\nu \in \mathbb{S}^{d-1}$.
    That this is in general not the case is shown in \cite{BenzoniGavage2007} where a particular system fails the continuous dependence of the eigenspaces with respect to $\nu$.
    However, for constantly hyperbolic systems, i.e. systems with constant multiplicity of the eigenvalues, the eigenspaces are continuous functions of $\nu \in \mathbb{S}^{d-1}$.
    Hence it is possible to find locally an eigenbasis of $\mathbf{A}^\ast(\nu)$, which depends continuously on $\nu$.
    \begin{rem}
        Here we additionally assume that the transformation is continuous with respect to each smooth boundary part and $\nu \in \mathbb{S}^{d-1}$.
        Following \cite{BenzoniGavage2007} for systems with constant coefficients the class of constantly hyperbolic systems, for example, enjoys this property.
    \end{rem}
    The boundary $\partial\Omega$ will be separated in the controllable and uncontrollable part, i.e. for $i = 1,\dots,n$
    \begin{align*}
        \Gamma_i^+ &:= \left\{\left.\mathbf{x} \in \partial\Omega\,\right|\,\lambda_i(\mathbf{x},\mathbf{n}(\mathbf{x})) \geq 0\right\},\\
        \Gamma_i^- &:= \left\{\left.\mathbf{x} \in \partial\Omega\,\right|\,\lambda_i(\mathbf{x},\mathbf{n}(\mathbf{x})) < 0\right\}.
    \end{align*}
    We further need to discuss a LMI associated to the system \eqref{eq:hyp_cons_sys} of the form
    \begin{align}
        \mathbf{A}(\mathbf{m}) := C\mathbf{Id} + \sum_{k=1}^dm_k\mathbf{A}^{(k)} \leq 0,\, C \in \R_{>0}\label{ineq:control_LMI}
    \end{align}
    which obviously may be reformulated to
    \begin{align}
        \overline{\mathbf{A}}(\mathbf{m}) := \mathbf{A}^{(0)} + \sum_{k=1}^d\overline{m}_k\mathbf{A}^{(k)} \geq 0, \label{ineq:control_LMI_v2}
    \end{align}
    with $\mathbf{A}^{(0)} = -\mathbf{Id}$ and where the coefficients are given by $\overline{m}_k = -m_k/C$. The matrices $\mathbf{A}^{(k)}$ are the Jacobians of system \eqref{eq:hyp_cons_sys}.
    In the case that the Jacobians of system \eqref{eq:hyp_cons_sys} are space dependent the coefficients $\mathbf{m}$ of the LMI may also depend on $\mathbf{x}$.
    Further it will turn out that for the boundary feedback control we need $\mathbf{m}$ to be the gradient of a function and hence we have the following result.
    \begin{lem}\label{defi:feas_potential}
        Assume there exists a function $\overline{\mu}(\mathbf{x})$ such that
        \begin{align}
            \overline{\mathbf{m}} := \nabla\overline{\mu}(\mathbf{x})\quad\text{and}\quad \overline{\mathbf{A}}(\mathbf{m}) := -\mathbf{Id} + \sum_{k=1}^d\overline{m}_k\mathbf{A}^{(k)} \geq 0.\label{cond:feas_potential}
        \end{align}
        Then system \eqref{eq:hyp_cons_sys} has a \emph{feasible Lyapunov potential} $\mu(\mathbf{x}) = -C\overline{\mu}(\mathbf{x})$ such that \eqref{ineq:control_LMI} holds.
    \end{lem}
    Of course the space dependent case is more intricate since we not only need any solution of the LMI, but also one which must be integrable.
    However, for constant matrices we have for such coefficients $\mathbf{m}$ that
    \[
      \mu(\mathbf{x}) = \sum_{k=1}^dm_kx_k
    \]
    is a suited potential.\\
    Finally we want to recall the definition of exponential stability, cf.\ \cite{Bastin2016,Hayat2021}.
    \begin{defi}
        A solution $\mathbf{w} \in C^1(0,T;H^s(\Omega))^n$ for $s\geq 1 + \frac{d}2$  of the initial boundary value problem \eqref{eq:hyp_cons_sys} is called \emph{exponentially stable} in the $L^2-$ sense, iff
        \begin{align}
            \|\mathbf{w}(t,.)\|_{L^2(\Omega)} \leq  C_1\exp(-C_2t)\|\mathbf{w}(0,.)\|_{L^2(\Omega)},\,C_1,C_2\in\R_{>0},\,t \in [0,+\infty).\label{def:asymp_stable}
        \end{align}
    \end{defi}
    The Lyapunov function we will introduce in Theorem \ref{thm:main} is equivalent to the $L^2(\Omega)$ norm.
    For later use, we define
    \begin{align*}
        \mathcal{A}(t,\mathbf{x}) &:=
        \left(\mathbf{w}^T(t,\mathbf{x})\mathbf{A}^{(1)}(\mathbf{x})\mathbf{w}(t,\mathbf{x}),\dots,
        \mathbf{w}^T(t,\mathbf{x})\mathbf{A}^{(d)}(\mathbf{x})\mathbf{w}(t,\mathbf{x})\right)^T\exp(\mu(\mathbf{x})),\\
        m_k(\mathbf{x}) &:= \p{}{x_k}\mu(\mathbf{x}),\,k=1,\dots,d.
    \end{align*}
    Now, we state the main result of this paper.
    \begin{thm}\label{thm:main}
        Let $\mathbf{w}(t,\mathbf{x}) \in C^1\left((0,T),H^s(\Omega)\right)^n$, $s \geq 1+ d/2,$ be a solution to the IBVP \eqref{eq:hyp_cons_sys}.
        Assume the problem \eqref{eq:hyp_cons_sys} admits a feasible Lyapunov potential.
        Thus there exist $\mu(\mathbf{x}) \in H^s(\Omega)$ such that
        \begin{align}
            C_A\mathbf{Id} + \sum_{k=1}^d m_k(\mathbf{x})\mathbf{A}^{(k)}(\mathbf{x}) \leq 0\label{eq:mu_pde}
        \end{align}
        holds true for some value $C_A \in \R_{>0}$ with $C_A > C_B$ and $C_B$ given as in \eqref{ineq:source_est}.
        Define the Lyapunov function as follows
        \begin{align}
            L(t) = \int_\Omega \mathbf{w}(t,\mathbf{x})^T\mathbf{w}(t,\mathbf{x})\exp(\mu(\mathbf{x}))\,\dd\mathbf{x}.\label{eq:lyapunov_gen}
        \end{align}
        Let the boundary condition for \eqref{eq:hyp_cons_sys} be given by
        \begin{align*}
            \mathbf{w}^{(i)}_{BC}(t,\mathbf{x}) = u_i(t,\mathbf{x}),\,\mathbf{x}\in\Gamma_i^-
        \end{align*}
        where the imposed feedback controls $\mathbf{u} = (u_1,\dots,u_n)^T$ satisfy $\mathbf{u} = \mathbf{T}\tilde{\mathbf{u}}$
        \begin{align}
            \begin{split}\label{ineq:boundary_control}
                &-\sum_{i=1}^n\int_{\Gamma_i^-}\lambda_i(\mathbf{x},\mathbf{n}(\mathbf{x}))\tilde{u}_i(t,\mathbf{x})^2\exp(\mu(\mathbf{x}))\,\dd\mathbf{x}
                \leq \sum_{i=1}^n\int_{\Gamma_i^+}\lambda_i(\mathbf{x},\mathbf{n}(\mathbf{x}))q_i(t,\mathbf{x})^2\exp(\mu(\mathbf{x}))\,\dd\mathbf{x},
            \end{split}
        \end{align}
        with $\mathbf{q} = \mathbf{T}^T\mathbf{w}$ and \eqref{eq:diag_transform}.\\
        Then, it holds for the Lyapunov function that for $0 < C_L = C_A - C_B$
        \begin{align*}
            \frac{\dd}{\dd t}L(t) \leq -C_LL(t).
        \end{align*}
        Further, the solution $\mathbf{w}$ of the IBVP \eqref{eq:hyp_cons_sys} is exponentially stable in the $L^2-$sense \eqref{def:asymp_stable}.
    \end{thm}
    The proof is given in the following section for better readability. We refer to Remark \ref{rem:control} for particular cases of interest and conditions under which $u_i$ exists.
    Note that the choice of $s$ in the space $H^s$ guarantees that $\mu$ and $\mathbf{w}$ are differentiable in $\mathbf{x}$ by Sobolev embedding.
    \section{Proof of Theorem \ref{thm:main}}\label{sec:prf}
    \begin{prf}
        We want to show that the Lyapunov function \eqref{eq:lyapunov_gen} decays exponentially fast under the given conditions.
        We obtain for the time derivative of the Lyapunov function
        \begin{align*}
            \frac{\dd}{\dd t}L(t) &= \frac{\dd}{\dd t}\int_\Omega \mathbf{w}^T(t,\mathbf{x})\mathbf{w}(t,\mathbf{x})\exp(\mu(\mathbf{x}))\,\dd\mathbf{x}
            = 2\int_\Omega \mathbf{w}^T(t,\mathbf{x})\p{}{t}\mathbf{w}(t,\mathbf{x})\exp(\mu(\mathbf{x}))\,\dd\mathbf{x}.
        \end{align*}
        Replacing the time derivative of $\mathbf{w}$ using the PDE system \eqref{eq:hyp_cons_sys} gives
        \begin{align*}
            \frac{\dd}{\dd t}L(t) &= 2\int_\Omega \mathbf{w}^T(t,\mathbf{x})
            \left[-\sum_{k=1}^d\mathbf{A}^{(k)}(\mathbf{x})\p{}{x_k}\mathbf{w}(t,\mathbf{x}) - \mathbf{B}(\mathbf{x})\mathbf{w}(t,\mathbf{x})\right]\exp(\mu(\mathbf{x}))\,\dd\mathbf{x}\\
            &= -2\int_\Omega \sum_{k=1}^d\mathbf{w}^T(t,\mathbf{x})\mathbf{A}^{(k)}(\mathbf{x})\p{}{x_k}\mathbf{w}(t,\mathbf{x})\exp(\mu(\mathbf{x}))
            + \mathbf{w}^T(t,\mathbf{x})\mathbf{B}(\mathbf{x})\mathbf{w}(t,\mathbf{x})\exp(\mu(\mathbf{x}))\,\dd\mathbf{x}
        \end{align*}
        Now with using the product rule and introducing the abbreviations $\mathcal{A}$ and $m_k(\mathbf{x})$ we obtain
        \begin{align*}
            &\frac{\dd}{\dd t}L(t) = \dots\\
            &= -\int_\Omega \sum_{k=1}^d\left[\p{}{x_k}\left(\mathbf{w}^T(t,\mathbf{x})\mathbf{A}^{(k)}(\mathbf{x})\mathbf{w}(t,\mathbf{x})\exp(\mu(\mathbf{x}))\right)
            - \mathbf{w}^T(t,\mathbf{x})\left(\p{}{x_k}\mathbf{A}^{(k)}(\mathbf{x})\right)\mathbf{w}(t,\mathbf{x})\exp(\mu(\mathbf{x}))\right.\\
            &-\left. m_k(\mathbf{x})\mathbf{w}^T(t,\mathbf{x})\mathbf{A}^{(k)}(\mathbf{x})\mathbf{w}(t,\mathbf{x})\exp(\mu(\mathbf{x}))\right]
            + 2\mathbf{w}^T(t,\mathbf{x})\mathbf{B}(\mathbf{x})\mathbf{w}(t,\mathbf{x})\exp(\mu(\mathbf{x}))\,\dd\mathbf{x}\\
            &= -\int_\Omega \nabla\cdot\mathcal{A}(\mathbf{x})\,\dd\mathbf{x}
            + \int_\Omega \mathbf{w}^T(t,\mathbf{x})\sum_{k=1}^d\left(m_k(\mathbf{x})\mathbf{A}^{(k)}(\mathbf{x}) + \p{}{x_k}\mathbf{A}^{(k)}(\mathbf{x})\right)\mathbf{w}(t,\mathbf{x})\exp(\mu(\mathbf{x}))\\
            &- 2\mathbf{w}^T(t,\mathbf{x})\mathbf{B}(\mathbf{x})\mathbf{w}(t,\mathbf{x})\exp(\mu(\mathbf{x}))\,\dd\mathbf{x}\\
            &= -\underbrace{\int_{\partial\Omega} \mathcal{A}(t,\mathbf{x})\cdot\mathbf{n}(\mathbf{x})\,\dd\mathbf{x}}_{=:\mathcal{B}(t)}\\
            &+ \underbrace{\int_\Omega \mathbf{w}^T(t,\mathbf{x})
            \left[\sum_{k=1}^d\left(m_k(\mathbf{x})\mathbf{A}^{(k)}(\mathbf{x}) + \p{}{x_k}\mathbf{A}^{(k)}(\mathbf{x})\right)- 2\mathbf{B}(\mathbf{x})\right]
            \mathbf{w}(t,\mathbf{x})\exp(\mu(\mathbf{x}))\,\dd\mathbf{x}}_{=: \mathcal{I}(t)}.
        \end{align*}
        In the following we split the proof into two parts. First we show that the boundary term $\mathcal{B}(t)$ is non-negative under the given conditions.
        In the second part we estimate $\mathcal{I}(t)$ such that we can use Gronwalls Lemma to show the decay of the Lyapunov function.
        \subsection{Estimating the Boundary Integral}\label{sec:BC}
        In the following we want to discuss the boundary term. We have
        \begin{align*}
            \mathcal{A}(t,\mathbf{x})\cdot\mathbf{n}(\mathbf{x}) &= \sum_{k=1}^d\mathbf{w}^T(t,\mathbf{x})\mathbf{A}^{(k)}\mathbf{w}(t,\mathbf{x})n_k(\mathbf{x})\exp(\mu(\mathbf{x}))\\
            &= \mathbf{w}^T(t,\mathbf{x})\left(\sum_{k=1}^dn_k(\mathbf{x})\mathbf{A}^{(k)}(\mathbf{x})\right)\mathbf{w}(t,\mathbf{x})\exp(\mu(\mathbf{x}))\\
            &= \mathbf{w}^T(\mathbf{x})\mathbf{A}^\ast(\mathbf{x},\mathbf{n}(\mathbf{x}))\mathbf{w}(t,\mathbf{x})\exp(\mu(\mathbf{x})).
        \end{align*}
        Since the system is symmetric we can transform the quadratic form according to \eqref{eq:diag_transform} and thus have
        \begin{align*}
            \mathcal{B}(t) = \int_{\partial\Omega} \mathcal{A}(t,\mathbf{x})\cdot\mathbf{n}(\mathbf{x})\,\dd\mathbf{x}
            = \int_{\partial\Omega}\mathbf{q}^T(t,\mathbf{x})\mathbf{\Lambda}^\ast(\mathbf{x},\mathbf{n}(\mathbf{x}))\mathbf{q}(t,\mathbf{x})\exp(\mu(\mathbf{x}))\,\dd\mathbf{x}.
        \end{align*}
        Using the partitioning of the boundary introduced before we rewrite the boundary term and obtain
        \begin{align*}
            \mathcal{B}(t) &= \int_{\partial\Omega} \mathcal{A}(\mathbf{x})\cdot\mathbf{n}(\mathbf{x})\,\dd\mathbf{x}
            = \int_{\partial\Omega}\mathbf{q}^T(t,\mathbf{x})\mathbf{\Lambda}^\ast(\mathbf{x},\mathbf{n}(\mathbf{x}))\mathbf{q}(t,\mathbf{x})\exp(\mu(\mathbf{x}))\,\dd\mathbf{x}\\
            &= \sum_{i=1}^n\int_{\Gamma_i^-}\lambda_i(\mathbf{x},\mathbf{n}(\mathbf{x}))\exp(\mu(\mathbf{x}))\tilde{u}_i(t,\mathbf{x})^2\,\dd\mathbf{x}
            + \sum_{i=1}^n\int_{\Gamma_i^+}\lambda_i(\mathbf{x},\mathbf{n}(\mathbf{x}))\exp(\mu(\mathbf{x}))q_i(t,\mathbf{x})^2\,\dd\mathbf{x}\\
            &\geq 0.
        \end{align*}
        Due to $\mathbf{q} = \mathbf{T}^T\mathbf{w}$ we have $\mathbf{u} = \mathbf{T}\tilde{\mathbf{u}}$.
        \subsection{Estimating the Volume Integral}
        We now want to study the expression
        \begin{align*}
            \mathcal{I}(t) &= \int_\Omega \mathbf{w}^T(t,\mathbf{x})
            \left[\sum_{k=1}^d\left(m_k(\mathbf{x})\mathbf{A}^{(k)}(\mathbf{x}) + \p{}{x_k}\mathbf{A}^{(k)}(\mathbf{x})\right)- 2\mathbf{B}(\mathbf{x})\right]
            \mathbf{w}(t,\mathbf{x})\exp(\mu(\mathbf{x}))\,\dd\mathbf{x}\\
            &= \int_\Omega \mathbf{w}^T(t,\mathbf{x})\left(\sum_{k=1}^dm_k(\mathbf{x})\mathbf{A}^{(k)}(\mathbf{x})\right)\mathbf{w}(t,\mathbf{x})\exp(\mu(\mathbf{x}))\,\dd\mathbf{x}\\
            &+ \int_\Omega \mathbf{w}^T(t,\mathbf{x})\left(\sum_{k=1}^d\p{}{x_k}\mathbf{A}^{(k)}(\mathbf{x})- 2\mathbf{B}(\mathbf{x})\right)\mathbf{w}(t,\mathbf{x})\exp(\mu(\mathbf{x}))\,\dd\mathbf{x}
        \end{align*}
        We first want to discuss the term
        \begin{align}
            \mathbf{w}^T(t,\mathbf{x})\left(\sum_{k=1}^d\p{}{x_k}\mathbf{A}^{(k)}(\mathbf{x}) - 2\mathbf{B}(\mathbf{x})\right)\mathbf{w}(t,\mathbf{x}).\label{reminder_term}
        \end{align}
        If \eqref{reminder_term} is less or equal to zero there is nothing to do and it can simply be estimated by zero.
        Assume \eqref{reminder_term} is strictly positive. In the case that $\mathbf{B}$ is not symmetric it can be replaced by its symmetric part $\mathbf{B}^{sym} = (\mathbf{B} + \mathbf{B}^T)/2$.
        Since the $\mathbf{A}^{(k)}$ are symmetric this property is inherited by the derivatives of these matrices.
        Hence the following estimate can be established
        \begin{align}
            \mathbf{w}^T(t,\mathbf{x})\left(\sum_{k=1}^d\p{}{x_k}\mathbf{A}^{(k)}(\mathbf{x}) - 2\mathbf{B}^{sym}(\mathbf{x})\right)\mathbf{w}(t,\mathbf{x}) \leq C_B\|\mathbf{w}(t,\mathbf{x})\|_2^2,\,C_B\in\R_{>0}.
            \label{ineq:source_est}
        \end{align}
        Furthermore, by \eqref{eq:mu_pde}, $\mu(\mathbf{x})$ and thus the $m_k(\mathbf{x})$, respectively, are chosen such that the linear matrix inequality
        \begin{align*}
            C_A\mathbf{Id} + \sum_{k=1}^d m_k(\mathbf{x})\mathbf{A}^{(k)}(\mathbf{x}) \leq 0
        \end{align*}
        holds.
        Thus we have for $C_A > C_B$ and $0 < C_L \leq C_A - C_B$
        \begin{align*}
            \mathcal{I}(t) &= \int_\Omega \mathbf{w}^T(t,\mathbf{x})
            \left[\sum_{k=1}^d\left(m_k(\mathbf{x})\mathbf{A}^{(k)}(\mathbf{x}) + \p{}{x_k}\mathbf{A}^{(k)}(\mathbf{x})\right)- 2\mathbf{B}(\mathbf{x})\right]
            \mathbf{w}(t,\mathbf{x})\exp(\mu(\mathbf{x}))\,\dd\mathbf{x}\\
            &\leq \int_\Omega \mathbf{w}^T(t,\mathbf{x})\left[C_B\mathbf{Id} + \sum_{k=1}^d m_k(\mathbf{x})\mathbf{A}^{(k)}(\mathbf{x})\right]\mathbf{w}(t,\mathbf{x})\exp(\mu(\mathbf{x}))\,\dd\mathbf{x}\\
            &\leq -C_L\int_\Omega \mathbf{w}^T(t,\mathbf{x})\mathbf{w}(t,\mathbf{x})\exp(\mu(\mathbf{x}))\,\dd\mathbf{x}.
        \end{align*}
        Altogether we yield for the Lyapunov function
        \begin{align*}
            \frac{\dd}{\dd t}L(t) = -\mathcal{B}(t) + \mathcal{I}(t) \leq \mathcal{I}(t) \leq -C_L\int_\Omega \mathbf{w}^T(t,\mathbf{x})\mathbf{w}(t,\mathbf{x})\exp(\mu(\mathbf{x}))\,\dd\mathbf{x} = -C_LL(t).
        \end{align*}
        Applying Gronwall's Lemma gives the claimed exponential decay and thus the proof is complete.
    \end{prf}
    \begin{rem}
        We want to remark that the obtained result can be extended to the exponential stability in the corresponding $H^s$-norm.
        One just has to add other appropriate terms in the definition of the Lyapunov function. The $H^s$-norm of a function $u \in H^s(\Omega)$ is given by
        \[
          \|u\|_{H^s(\Omega)} := \left(\sum_{|\alpha|\leq s}\|\partial^\alpha u\|^2_2\right)^{\frac{1}{2}}.
        \]
        For vector-valued functions the norm is interpreted component wise.
        Now we define the Lyapunov function accordingly
        \begin{align}
            L(t) := \sum_{|\alpha|\leq s}\left\|\exp\left(\frac{1}{2}\mu^{(\alpha)}\right)\partial^\alpha \mathbf{u}\right\|^2_2.\label{def:Hs_Lyapunov}
        \end{align}
        Note that by introducing the $\exp\left(\mu^{(\alpha)}\right)$ we want to highlight that in principle it is possible to choose suited Lyapunov potentials for every occurring multi-index $|\alpha| \leq s$.
        Each Lyapunov potential has to satisfy a LMI similar to \eqref{ineq:control_LMI} with respect to the system corresponding to the studied multi-index.
        It has to be studied under which conditions this is possible for a given system.
        Once the Lyapunov function is introduced the strategy is the same and due to the similar structure inherited by the governing PDE \eqref{eq:hyp_cons_sys} the estimates can be adopted from the previous proof.
        We exemplary want to give the Lyapunov function in the case $H^1$
        \begin{align*}
            L(t) &:= \int_\Omega \mathbf{w}(t,\mathbf{x})^T\mathbf{w}(t,\mathbf{x})\exp\left(\mu^{(0)}(\mathbf{x})\right)\,\dd\mathbf{x}\\
            &+ \sum_{k=1}^d\int_\Omega \left(\partial_{x_k}\mathbf{w}(t,\mathbf{x})\right)^T\left(\partial_{x_k}\mathbf{w}(t,\mathbf{x})\right)\exp\left(\mu^{(k)}(\mathbf{x})\right)\,\dd\mathbf{x}
        \end{align*}
    \end{rem}
    \begin{rem}\label{rem:lmi}
        During the proof we encounter the term \eqref{reminder_term}.
        In particular in the proof we only treat the situation when this expression is strictly positive since under the given assumptions, namely \eqref{eq:mu_pde},
        it can simply be estimated by zero in the case that it is less or equal to zero.
        However, there may be situations when the estimate of the volume integral may benefit from this term.
        In the case that the reminder \eqref{reminder_term} should be taken into account we introduce
        \[
          \mathcal{R}(\mathbf{x}) := \sum_{k=1}^d\p{}{x_k}\mathbf{A}^{(k)}(\mathbf{x}) - 2\mathbf{B}^{sym}(\mathbf{x})
        \]
        and demand
        \begin{align}
            \mathbf{A}(\mathbf{m}) := C\mathbf{Id} + \mathcal{R}(\mathbf{x}) + \sum_{k=1}^dm_k(\mathbf{x})\mathbf{A}^{(k)}(\mathbf{x}) \leq 0,\, C \in \R_{>0}.\label{ineq:control_LMI2}
        \end{align}
        The LMI \eqref{ineq:control_LMI2} then replaces \eqref{eq:mu_pde}.
    \end{rem}
    \begin{rem}\label{rem:control}
        The inequality for the (transformed) control $\tilde{u}_i$ can be simplified under additional assumptions.
        \begin{enumerate}[(i)]
            \item A possible simplifying assumption on the controls is that they are uniform in space, i.e.\ $\tilde{u}_i = \tilde{u}_i(t)$.
            This leads to the following condition for the controls
            \begin{align*}
                -\sum_{i=1}^n\tilde{u}_i(t)^2\int_{\Gamma_i^-}\lambda_i(\mathbf{x},\mathbf{n}(\mathbf{x}))e^{\mu(\mathbf{x})}\,\dd\mathbf{x}
                \leq \sum_{i=1}^n\int_{\Gamma_i^+}q_i(t,\mathbf{x})^2\lambda_i(\mathbf{x},\mathbf{n}(\mathbf{x}))e^{\mu(\mathbf{x})}\,\dd\mathbf{x}.
            \end{align*}
            \item Further the same control may be applied to all components which leads to $\tilde{u}_i \equiv \tilde{u}(t,\mathbf{x}),\,i=1,\dots,n$.
            This leads to the following condition for the controls
            \begin{align*}
                -\sum_{i=1}^n\int_{\Gamma^-_i}\tilde{u}(t,\mathbf{x})^2\lambda_i(\mathbf{x},\mathbf{n}(\mathbf{x}))e^{\mu(\mathbf{x})}\,\dd\mathbf{x}
                \leq \sum_{i=1}^n\int_{\Gamma_i^+}q_i(t,\mathbf{x})^2\lambda_i(\mathbf{x},\mathbf{n}(\mathbf{x}))e^{\mu(\mathbf{x})}\,\dd\mathbf{x}.
            \end{align*}
            \item Assume that both previous assumptions hold, i.e. $\tilde{u}_i \equiv \tilde{u}(t),\,i=1,\dots,n$. Then we can give an explicit condition for the feedback control as
            \begin{align}
                \tilde{u}(t)^2 \leq -\left(\sum_{i=1}^n\int_{\Gamma^-_i}\lambda_i(\mathbf{x},\mathbf{n}(\mathbf{x}))e^{\mu(\mathbf{x})}\,\dd\mathbf{x}\right)^{-1}
                \sum_{i=1}^n\int_{\Gamma_i^+}q_i(t,\mathbf{x})^2\lambda_i(\mathbf{x},\mathbf{n}(\mathbf{x}))e^{\mu(\mathbf{x})}\,\dd\mathbf{x}.\label{ineq:boundary_control2}
            \end{align}
            Thus a possible control would be
            \begin{align}
                \tilde{u}(t) = C\sqrt{-\left(\sum_{i=1}^n\int_{\Gamma^-_i}\lambda_i(\mathbf{x},\mathbf{n}(\mathbf{x}))e^{\mu(\mathbf{x})}\,\dd\mathbf{x}\right)^{-1}
                \sum_{i=1}^n\int_{\Gamma_i^+}q_i(t,\mathbf{x})^2\lambda_i(\mathbf{x},\mathbf{n}(\mathbf{x}))e^{\mu(\mathbf{x})}\,\dd\mathbf{x}}\label{eq:boundary_control}
            \end{align}
            for $C\in[-1,1]$.
        \end{enumerate}
    \end{rem}
    \section{Application to Hyperbolic Systems}\label{sec:app_hyp_sys}
    In this section we discuss the application of the stabilization result.
    As a first example we want to mention that recent results have shown the applicability of this approach to Hamilton-Jacobi equations, see \cite{Bambach2022,Herty2023}.
    The study on the stabilization of Hamilton--Jacobi equation has been motivated by a stabilization problem arising in a forming process and we refer to the afore mentioned literature for further details.\\
    We now want to study the general case and show that the linearization matches \eqref{eq:hyp_cons_sys}. Let us consider the following hyperbolic system
    \begin{align}
        \begin{dcases}
            \p{}{t}\mathbf{U}(t,\mathbf{x}) + \mathbf{F}(\mathbf{x},\mathbf{U}(t,\mathbf{x}),\mathbf{U}_{x_1}(t,\mathbf{x}),\dots,\mathbf{U}_{x_d}(t,\mathbf{x})) &= 0,\;(t,\mathbf{x})\in(0,T)\times\Omega,\\
            \mathbf{U}(0,\mathbf{x}) &= \mathbf{U}_0(\mathbf{x}),\;\mathbf{x}\in\Omega,\\
            \mathbf{U}(t,\mathbf{x}) &= \mathbf{U}_{BC}(t,\mathbf{x}),\;(t,\mathbf{x})\in[0,T)\times\partial\Omega
        \end{dcases}\label{eq:gen_hyp_pde}
    \end{align}
    Here we have introduced the function
    \begin{align*}
        &\mathbf{F}: \R^d\times\underbrace{\R^n\times\cdots\times\R^n}_{d+1} \to \R^n\\
        &(\mathbf{x},\mathbf{u},\mathbf{p}_1,\dots,\mathbf{p}_d) \mapsto \mathbf{F}(\mathbf{x},\mathbf{u},\mathbf{p}_1,\dots,\mathbf{p}_d) \equiv \mathbf{F}(\mathbf{x},\mathbf{u},\mathbf{p}).
    \end{align*}
    We want to linearize system \eqref{eq:gen_hyp_pde} near a steady state solution $\overline{\mathbf{U}}(\mathbf{x})$ of \eqref{eq:gen_hyp_pde}.
    Thus we have to linearize $\mathbf{F}$ near a reference state $(\overline{\mathbf{u}},\overline{\mathbf{p}})\equiv(\overline{\mathbf{u}},\overline{\mathbf{p}}_1,\dots,\overline{\mathbf{p}}_d)$
    \begin{align}
        \mathbf{F}(\mathbf{x},\mathbf{u},\mathbf{p}) &\approx \mathbf{F}(\mathbf{x},\overline{\mathbf{u}},\overline{\mathbf{p}}) +
        \underbrace{\p{}{\mathbf{u}}\mathbf{F}(\mathbf{x},\overline{\mathbf{u}},\overline{\mathbf{p}})}_{=:\mathbf{B}(\mathbf{x})}(\mathbf{u} - \overline{\mathbf{u}}) +
        \p{}{\mathbf{p}}\mathbf{F}(\mathbf{x},\overline{\mathbf{u}},\overline{\mathbf{p}})(\mathbf{p} - \overline{\mathbf{p}})\notag\\
        &= \mathbf{F}(\mathbf{x},\overline{\mathbf{u}},\overline{\mathbf{p}}) + \mathbf{B}(\mathbf{x})(\mathbf{u} - \overline{\mathbf{u}}) +
        \sum_{k=1}^d\underbrace{\p{}{\mathbf{p}_k}\mathbf{F}(\mathbf{x},\overline{\mathbf{u}},\overline{\mathbf{p}})}_{=:\mathbf{A}^{(k)}(\mathbf{x})}(\mathbf{p}_k - \overline{\mathbf{p}}_k)\notag\\
        &=\mathbf{F}(\mathbf{x},\overline{\mathbf{u}},\overline{\mathbf{p}}) + \mathbf{B}(\mathbf{x})(\mathbf{u} - \overline{\mathbf{u}}) +
        \sum_{k=1}^d\mathbf{A}^{(k)}(\mathbf{x})(\mathbf{p}_k - \overline{\mathbf{p}}_k).\label{approx:F}
    \end{align}
    We now consider deviations of this desired state
    \[
      \mathbf{U}(t,\mathbf{x})  = \mathbf{\overline{U}}(\mathbf{x}) + \mathbf{w}(t,\mathbf{x})
    \]
    and for some initial disturbance we assume $\mathbf{w}(0,\mathbf{x}) = (w_1(0,\mathbf{x}),\dots,w_n(0,\mathbf{x})) \neq 0$.
    By the following computation we obtain a linear system using \eqref{approx:F} starting from \eqref{eq:gen_hyp_pde}
    \begin{align*}
        0 &= \p{}{t}\mathbf{U}(t,\mathbf{x}) + \mathbf{F}(\mathbf{x},\mathbf{U}(t,\mathbf{x}),\mathbf{U}_{x_1}(t,\mathbf{x}),\dots,\mathbf{U}_{x_d}(t,\mathbf{x}))\\
        &\approx \p{}{t}\mathbf{U}(t,\mathbf{x}) + \mathbf{F}(\mathbf{x},\overline{\mathbf{U}}(\mathbf{x}),\overline{\mathbf{U}}_{x_1}(\mathbf{x}),\dots,\overline{\mathbf{U}}_{x_d}(\mathbf{x}))
        + \mathbf{B}(\mathbf{x})(\mathbf{U}(t,\mathbf{x}) - \overline{\mathbf{U}}(\mathbf{x}))\\ &+ \sum_{k=1}^d\mathbf{A}^{(k)}(\mathbf{x})(\mathbf{U}_{x_k}(t,\mathbf{x}) - \overline{\mathbf{U}}_{x_k}(\mathbf{x}))\\
        &= \p{}{t}(\mathbf{\overline{U}}(\mathbf{x}) + \mathbf{w}(t,\mathbf{x}))
        + \mathbf{F}(\mathbf{x},\overline{\mathbf{U}}(\mathbf{x}),\overline{\mathbf{U}}_{x_1}(\mathbf{x}),\dots,\overline{\mathbf{U}}_{x_d}(\mathbf{x}))
        + \mathbf{B}(\mathbf{x})\mathbf{w}(t,\mathbf{x})\\
        &+ \sum_{k=1}^d\mathbf{A}^{(k)}(\mathbf{x})\p{}{x_k}\mathbf{w}(t,\mathbf{x})\\
        &= \p{}{t}\mathbf{w}(t,\mathbf{x}) + \sum_{k=1}^d\mathbf{A}^{(k)}(\mathbf{x})\p{}{x_k}\mathbf{w}(t,\mathbf{x}) + \mathbf{B}(\mathbf{x})\mathbf{w}(t,\mathbf{x})
    \end{align*}
    Thus the linearized system has the desired form to apply the above results, highlighting the wide range of applicability of the present work.
    \subsection{Barotropic Euler Equations}\label{sec:halfspace_prob}
    The barotropic Euler system (isothermal or isentropic) is of huge interest in the literature. The system reads in two dimensions
    \begin{align}
        \begin{split}\label{baro_eulereq_sys}
            \frac{\partial}{\partial t}\rho + \nabla_x\cdot\left(\rho\mathbf{v}\right) &= 0,\\
            \frac{\partial}{\partial t}(\rho v_i) + \nabla_x\cdot\left(\rho v_i\mathbf{v} + p\mathbf{e}^{(i)}\right) &= 0,\;i \in \{1,2\}
        \end{split}
    \end{align}
    where $\rho > 0$ denotes the mass density, $\mathbf{v} = (v_1,v_2)^T$ is the velocity and $p$ is the pressure as function of the density.
    The fluxes are given by
    \begin{align}
        \mathbf{F}^{(1)}(\mathbf{U}) = \begin{pmatrix} \rho v_1\\ \rho v_1^2 + p\\ \rho v_1v_2\end{pmatrix},\quad
        \mathbf{F}^{(2)}(\mathbf{U}) = \begin{pmatrix} \rho v_2\\ \rho v_1v_2\\ \rho v_2^2 + p\end{pmatrix}.
        \label{baro_eulereq:cons_flux}
    \end{align}
    Thus we can write the system in the following conservative form
    \[
      \p{}{t}\mathbf{U} + \p{}{x_1}\mathbf{F}^{(1)}(\mathbf{U}) + \p{}{x_2}\mathbf{F}^{(2)}(\mathbf{U}) = 0
    \]
    and hence obtain the following linearization
    \[
      \p{}{t}\mathbf{w} + \mathbf{A}^{(1)}(\overline{\mathbf{U}})\p{}{x_1}\mathbf{w} + \mathbf{A}^{(2)}(\overline{\mathbf{U}})\p{}{x_2}\mathbf{w} = 0.
    \]
    Now we could apply the obtained results to the linearized conservative system. However, since we assume some regularity of the solution and the deviation, it is beneficial to use a suited set of primitive variables.
    Here we will use $(\rho,\mathbf{v})$ which results in the following systems, cf.\ \cite{Warnecke1999,Toro2009},
    \begin{align}
        \begin{split}\label{baro_eulereq_prim_sys}
            \frac{\partial}{\partial t}\rho + \mathbf{v}\cdot\nabla_x\rho + \rho\nabla_x\cdot\mathbf{v} &= 0,\\
            \frac{\partial}{\partial t}v_i + \mathbf{v}\cdot\nabla_xv_i  + \frac{a^2}{\rho}\p{}{x_i}\rho &= 0,\;i \in \{1,2\}
        \end{split}
    \end{align}
    with $a > 0$ denoting the speed of sound.
    The Jacobians are given by
    \begin{align}
        \mathbf{A}^{(1)} = \begin{pmatrix}v_1 & \rho & 0\\ \frac{a^2}{\rho} & v_1 & 0\\ 0 & 0 & v_1\end{pmatrix}\quad\text{and}\quad
        \mathbf{A}^{(2)} = \begin{pmatrix}v_2 & 0 & \rho\\ 0 & v_2 & 0\\ \frac{a^2}{\rho} & 0 & v_2\end{pmatrix}.\label{baro_eulereq_jacobians_prim}
    \end{align}
    When linearizing the system we simply fix the Jacobians at the desired reference state.
    We consider the reference state $\overline{\rho} \in \R_{>0}, \overline{\mathbf{v}}\in\R^2\setminus\{\mathbf{0}\}$.
    The linearized system can be symmetrized by applying the following variable transformation $\mathbf{w}\to\mathbf{w}^\ast=(r,\mathbf{v})^T$
    \[
      \rho \to r := \frac{\overline{a}}{\overline{\rho}}\rho,
    \]
    see also \cite{Godlewski1996}. We hence obtain
     \[
      \p{}{t}\mathbf{w}^\ast + \overline{\mathbf{A}}^{(1)}\p{}{x_1}\mathbf{w}^\ast + \overline{\mathbf{A}}^{(2)}\p{}{x_2}\mathbf{w}^\ast = 0
    \]
    with the Jacobians
    \begin{align}
        \overline{\mathbf{A}}^{(1)} = \begin{pmatrix}\overline{v}_1 & \overline{a} & 0\\ \overline{a} & \overline{v}_1 & 0\\ 0 & 0 & \overline{v}_1\end{pmatrix}\quad\text{and}\quad
        \overline{\mathbf{A}}^{(2)} = \begin{pmatrix}\overline{v}_2 & 0 & \overline{a}\\ 0 & \overline{v}_2 & 0\\ \overline{a} & 0 & \overline{v}_2\end{pmatrix}.
        \label{baro_eulereq_jacobians_prim_diag_symm}
    \end{align}
    Since the Jacobians are constant and $\mathbf{B} \equiv \mathbf{0}$ we do not have to discuss the term \eqref{ineq:source_est}. Thus we have $C_B = 0$ and set $C_L := C_A > 0$.
    We now discuss the LMI
    \begin{align}
        C_L\mathbf{Id} + m_1\overline{\mathbf{A}}^{(1)} + m_2\overline{\mathbf{A}}^{(2)} \leq 0\label{ineq:LMI_euler}
    \end{align}
    which, as stated before, is equivalent to
    \begin{align*}
        \mathbf{A}(\overline{\mathbf{m}}) &:= -\mathbf{Id} + \overline{m}_1\overline{\mathbf{A}}^{(1)} + \overline{m}_2\overline{\mathbf{A}}^{(2)} \geq 0,\\
        \mathbf{A}(\overline{\mathbf{m}}) &=
        \begin{pmatrix}
            -1 + \overline{m}_1\overline{v}_1 + \overline{m}_2\overline{v}_2    & \overline{m}_1\overline{a}                                        & \overline{m}_2\overline{a}\\
            \overline{m}_1\overline{a}                                          & -1 + \overline{m}_1\overline{v}_1 + \overline{m}_2\overline{v}_2  & 0\\
            \overline{m}_2\overline{a}                                          & 0                                                                 & -1 + \overline{m}_1\overline{v}_1 + \overline{m}_2\overline{v}_2
        \end{pmatrix}.
    \end{align*}
    A sufficient criterion would be to choose $\overline{m}_1$ and $\overline{m}_2$ such that $\mathbf{A}(\overline{\mathbf{m}})$ becomes strongly diagonal dominant, i.e., the diagonal entries are strict positive and
    larger than the sum of the absolute row entries. Due to the special structure of the matrix we hence look for coefficients with the property that
    \[
      -1 + \overline{m}_1\overline{v}_1 + \overline{m}_2\overline{v}_2 = \sigma\left(|\overline{m}_1| + |\overline{m}_2|\right)\overline{a},\;\sigma\in(1,\infty).
    \]
    With this choice we have
    \[
      -1 + \overline{m}_1\overline{v}_1 + \overline{m}_2\overline{v}_2 > \left(|\overline{m}_1| + |\overline{m}_2|\right)\overline{a} \geq \max\{|\overline{m}_1|,|\overline{m}_2|\}\overline{a}
      \geq \min\{|\overline{m}_1|,|\overline{m}_2|\}\overline{a} \geq 0
    \]
    and hence the matrix $\mathbf{A}(\overline{\mathbf{m}})$ is positive definite.
    Therefore $-\mathbf{A}(\overline{\mathbf{m}})$ is negative definite and with the choice $\mathbf{m} = -C_L\overline{\mathbf{m}}$ the LMI \eqref{ineq:LMI_euler} holds and we further set
    \[
      \mu(\mathbf{x}) = m_1x_1 + m_2x_2.
    \]
    It remains to study the boundary term and prescribe boundary conditions, such that $\mathcal{BC} \geq 0$.
    To this end we will make use of the eigenstructure of the system given in the appendix for the readers convenience.
    For the sake of demonstration we consider the domain $\Omega = [0,1]\times[0,1]$ and a flow with $\overline{\mathbf{v}} = (\overline{v}_1,0)^T$ and $0 < \overline{a} < \overline{v}_1$.
    For the boundary integral we have
    \begin{align*}
        \mathcal{BC} &= \int_{\partial\Omega}\mathbf{w}^T\mathbf{A}^\ast(\mathbf{n})\mathbf{w}\exp(\mu(x_1,x_2))\,\dd\mathbf{x}
        = \int_{\partial\Omega}\mathbf{q}^T\mathbf{\Lambda}^\ast(\mathbf{n})\mathbf{q}\exp(\mu(x_1,y_2))\,\dd\mathbf{x}\\
        &= \sum_{i=1}^3\int_{\partial\Omega}q_i^2\lambda_i(\mathbf{n})\exp(\mu(x_1,x_2))\,\dd\mathbf{x}
    \end{align*}
    where $\mathbf{q} = \mathbf{T}^T(\mathbf{n})\mathbf{w}$ is calculated according to \eqref{euler_vec_trafo} and $\mathbf{n}$ denotes the outward pointing normal of the boundary $\partial\Omega$.
    Now we identify the controllable and uncontrollable parts of the boundary, i.e. for $i = 1,2,3$
    \begin{align*}
        \Gamma_i^+ &:= \left\{\left.\mathbf{x} \in \partial\Omega\,\right|\,\lambda_i(\mathbf{x},\mathbf{n}(\mathbf{x})) \geq 0\right\},\\
        \Gamma_i^- &:= \left\{\left.\mathbf{x} \in \partial\Omega\,\right|\,\lambda_i(\mathbf{x},\mathbf{n}(\mathbf{x})) < 0\right\}.
    \end{align*}
    These are given by
    \begin{align}
        \Gamma_1^- &= \partial\Omega\setminus\Gamma_1^+,\quad &\Gamma_1^+ &= \{1\}\times(0,1),\notag\\
        \Gamma_2^- &= \{0\}\times(0,1),\quad &\Gamma_2^+ &= \partial\Omega\setminus\Gamma_2^-,\label{bound_split}\\
        \Gamma_3^- &= \emptyset,\quad &\Gamma_3^+ &= \partial\Omega.\notag
    \end{align}
    Let us denote the general controls by $\varphi(t,x_1,x_2)$ for the first component and $\psi(t,x_2)$ for the second component, respectively.
    These have to be chosen such that
    \begin{align*}
        \mathcal{BC} &= \underbrace{\int_{\partial\Omega\setminus\Gamma_1^+}\varphi(t,x_1,x_2)^2(\mathbf{n}(x_1,x_2)\overline{v}_1 - \overline{a})\exp(\mu(x_1,x_2))\,\dd\mathbf{x} -
        \overline{v}_1\int_0^1\psi(t,x_2)^2\exp(\mu(0,x_2))\,\dd x_2}_{\leq 0}\\
        &+ \underbrace{(\overline{v}_1 - \overline{a})\int_0^1q_1^2\exp(\mu(1,x_2))\,\dd x_2 +
        \overline{v}_1\int_{\partial\Omega\setminus\Gamma_2^-}q_2^2\mathbf{n}(x_1,x_2)\exp(\mu(x_1,x_2))\,\dd\mathbf{x}}_{\geq 0}\\
        &+ \underbrace{\int_{\partial\Omega}q_3^2(\mathbf{n}(x_1,x_2)\overline{v}_1 + \overline{a})\exp(\mu(x_1,x_2))\,\dd\mathbf{x}}_{\geq 0} \geq 0.
    \end{align*}
    Thus one sufficient ansatz would be to choose the controls according to
    \begin{align*}
        &\int_{\partial\Omega\setminus\Gamma_1^+}\varphi(t,x_1,x_2)^2(\overline{a} - \mathbf{n}(x_1,x_2)\overline{v}_1)\exp(\mu(x_1,x_2))\,\dd\mathbf{x} \leq\textcolor{red}{\dots}\\
        &(\overline{v}_1 - \overline{a})\int_0^1q_1^2\exp(\mu(1,x_2))\,\dd x_2 + \int_{\partial\Omega}q_3^2(\mathbf{n}(x_1,x_2)\overline{v}_1 +
        \overline{a}\exp(\mu(x_1,x_2))\,\dd\mathbf{x}\\
        \text{and}\quad&\int_0^1\psi(t,x_2)^2\exp(\mu(0,x_2))\,\dd x_2 \leq \int_{\partial\Omega\setminus\Gamma_2^-}q_2^2\mathbf{n}(x_1,x_2)\exp(\mu(x_1,x_2))\,\dd\mathbf{x}.
    \end{align*}
    \section{Summary}
    A novel Lyapunov function for $L^2-$control of multi--dimensional systems of hyperbolic equations has been presented.
    We want to emphasize that up to our best knowledge this is the first result on stabilizing boundary feedback control applied to general multi--dimensional hyperbolic systems.
    The analysis relies on the fact that the system can be diagonalized on the boundary  and based on the proper Lyapunov function an estimate involving the symmetric Jacobians can be established.
    A stabilizing feedback control has been derived and exponential decay of a weighted $L^2-$norm has been proven.
    Illustrating examples have been presented to show the applicability of the presented result.
    Future work concerns numerical studies, the application to complex domains, the investigation and classification of problems which possess a feasible Lyapunov potential.
    It is of further interest to study whether this approach can also be generalized to the semilinear or the quasilinear case.
    \subsection*{Acknowledgments}
    \small{This research is part of the DFG SPP 2183 \emph{Eigenschaftsgeregelte Umformprozesse}, project 424334423.\\
    M.H. thanks the Deutsche Forschungsgemeinschaft (DFG, German Research Foundation)
    for the financial support through 525842915/SPP2410, 525853336/SPP2410, 320021702/GRK2326, 333849990/IRTG-2379, CRC1481, 423615040/SPP1962, 462234017, 461365406, ERS SFDdM035
    and under Germany’s Excellence Strategy EXC-2023 Internet of Production 390621612 and under the Excellence Strategy of the Federal Government and the L\"ander.\\
    F.T. thanks the Deutsche Forschungsgemeinschaft (DFG, German Research Foundation) for the financial support through 525939417/SPP2410.\\
    The authors gratefully acknowledge the valuable comments made by the referee.}
    \appendix
    \section{Eigenstructure of the barotropic Euler Equations}
    Subsequently we provide the detailed calculation for the eigenstructure of the considered system \eqref{baro_eulereq_jacobians_prim_diag_symm}.
    The system matrix is given by
    \begin{align}
        \mathbf{A}^\ast(\nu) = \nu_1\mathbf{A}^{(1)} + \nu_2\mathbf{A}^{(2)},\;\nu\in\mathbb{S}^2\label{bar_euler_sys_mat}
    \end{align}
    and thus we have
    \[
      \mathbf{A}^\ast(\nu) = \begin{pmatrix} \nu_1\overline{v}_1 + \nu_2\overline{v}_2  & \nu_1\overline{a}                         & \nu_2\overline{a}\\
                                             \nu_1\overline{a}                          & \nu_1\overline{v}_1 + \nu_2\overline{v}_2 & 0\\
                                             \nu_2\overline{a}                          & 0                                         & \nu_1\overline{v}_1 + \nu_2\overline{v}_2 \end{pmatrix}.
    \]
    The characteristic polynomial is given by
    \begin{align*}
        \chi(\lambda) &= (\lambda - (\nu_1\overline{v}_1 + \nu_2\overline{v}_2))^3 - (\lambda - (\nu_1\overline{v}_1 + \nu_2\overline{v}_2))\nu_1^2\overline{a}^2
        - (\lambda - (\nu_1\overline{v}_1 + \nu_2\overline{v}_2))\nu_2^2\overline{a}^2\\
        &= (\lambda - (\nu_1\overline{v}_1 + \nu_2\overline{v}_2))\left((\lambda - (\nu_1\overline{v}_1 + \nu_2\overline{v}_2))^2 - \left(\nu_1^2 + \nu_2^2\right)\overline{a}^2\right)\\
        &= (\lambda - (\nu_1\overline{v}_1 + \nu_2\overline{v}_2))\left((\lambda - (\nu_1\overline{v}_1 + \nu_2\overline{v}_2))^2 - \overline{a}^2\right)
    \end{align*}
    and hence we yield the following eigenvalues
    \begin{align}
        \lambda_1(\nu) &= \nu_1\overline{v}_1 + \nu_2\overline{v}_2 - \overline{a},\;\lambda_2(\nu) = \nu_1\overline{v}_1 + \nu_2\overline{v}_2,\;
        \lambda_3(\nu) = \nu_1\overline{v}_1 + \nu_2\overline{v}_2 + \overline{a}\label{eigvals_euler_sys}\\
        \text{with}\quad\lambda_1(\nu) &< \lambda_2(\nu) < \lambda_3(\nu).\notag
    \end{align}
    The corresponding right eigenvectors are thus obtained to be
    \begin{align}
        \mathbf{R}_1(\nu) = \frac{1}{\sqrt{2}}\begin{pmatrix} 1\\ -\nu_1\\ -\nu_2 \end{pmatrix},\quad
        \mathbf{R}_2(\nu) = \begin{pmatrix} 0\\ -\nu_2\\ \hphantom{-}\nu_1 \end{pmatrix}\quad\text{and}\quad
        \mathbf{R}_3(\nu) = \frac{1}{\sqrt{2}}\begin{pmatrix} 1\\ \nu_1\\ \nu_2 \end{pmatrix}.\label{eigvecs_euler_sys}
    \end{align}
    We hence have the following transformation matrix
    \[
      \mathbf{T}(\nu) = \frac{1}{\sqrt{2}}\begin{pmatrix} 1         & 0                         & 1\\
                                                          -\nu_1    & -\sqrt{2}\nu_2            & \nu_1\\
                                                          -\nu_2    & \hphantom{-}\sqrt{2}\nu_1 & \nu_2 \end{pmatrix}
    \]
    with $\mathbf{\Lambda}^\ast(\nu) = \mathbf{T}^T(\nu)\mathbf{A}^\ast(\nu)\mathbf{T}(\nu)$. The transformation of the state vector $\mathbf{w} = (r,v_1,v_2)^T$ is given by
    \begin{align}
        \mathbf{q} = \mathbf{T}^T(\nu)\mathbf{w} = \frac{1}{\sqrt{2}}\begin{pmatrix} r - (\nu_1v_1 + \nu_2v_2)\\ -\sqrt{2}(\nu_2v_1 - \nu_1v_2)\\ r + \nu_1v_1 + \nu_2v_2\end{pmatrix}
        \quad\text{and}\quad
        \mathbf{w} = \mathbf{T}(\nu)\mathbf{q} = \frac{1}{\sqrt{2}}\begin{pmatrix} q_1 + q_3\\ -\nu_1(q_1 - q_3) - \sqrt{2}\nu_2 q_2\\ -\nu_2(q_1 - q_3) + \sqrt{2}\nu_1 q_2 \end{pmatrix}
        \label{euler_vec_trafo}
    \end{align}
    %
    \phantomsection
    \bibliographystyle{abbrv}
    \bibliography{stb_hypsys_literature}
\end{document}